\documentclass[11pt]{article}

\usepackage[a4paper,margin=2.6cm]{geometry}
\usepackage{amsmath,amssymb,amsthm}
\usepackage{booktabs}
\usepackage{xcolor}
\usepackage{listings}
\usepackage[colorlinks=true,linkcolor=blue!60!black,citecolor=blue!60!black,urlcolor=blue!60!black]{hyperref}
\usepackage{microtype}
\usepackage{wasysym}
\usepackage{tikz}
\usetikzlibrary{arrows.meta,decorations.pathreplacing}
\tikzset{
  glyph/.style={line width=0.9pt, line cap=round, line join=round},
  pics/mpoint/.style={code={\fill (0,0) circle (0.06);}},
  pics/mline/.style={code={%
    \draw[glyph] (0,-0.5)--(0,0.5) (-0.28,0.12)--(0.28,0.12);}},
  pics/mcircle/.style={code={%
    \draw[glyph] (0,0) circle (0.36); \fill (0,0) circle (0.045);}},
  pics/mcrescent/.style={code={%
    \draw[glyph] (-0.4,0.2) arc (180:360:0.4);}},
  pics/maries/.style={code={%
    \draw[glyph] (0,-0.1) arc (0:180:0.2) (0,-0.1) arc (180:0:0.2);}},
  pics/mtaurus/.style={code={%
    \draw[glyph] (0,-0.2) circle (0.3) (-0.3,0.4) arc (180:360:0.3);}},
}

\hypersetup{
  pdftitle={Formalising John Dee's Monas Hieroglyphica in Lean 4},
  pdfauthor={Ivan Matychyn}}

\theoremstyle{definition}
\newtheorem{remark}{Remark}[section]

\definecolor{kwblue}{RGB}{0,0,180}
\definecolor{commentgreen}{RGB}{0,110,0}
\definecolor{stringmauve}{RGB}{140,20,140}
\definecolor{identred}{RGB}{150,0,0}

\lstdefinelanguage{Lean4}{
  keywords={theorem,lemma,def,abbrev,namespace,end,inductive,structure,
    private,deriving,by,intro,exact,refine,have,rw,cases,induction,simp,
    omega,decide,ac\_rfl,rfl,fun,match,with,let,in,at,
    example,package,import,open},
  keywordstyle=\bfseries\color{kwblue},
  sensitive=false,
  morecomment=[l]{--},
  morecomment=[l]{/-/},
  morecomment=[s]{/-}{-/},
  morestring=[b]",
  showstringspaces=false,
}
\lstdefinestyle{lean}{
  language=Lean4,
  basicstyle=\ttfamily\small,
  commentstyle=\itshape\color{commentgreen},
  stringstyle=\color{stringmauve},
  numbers=left,
  numberstyle=\tiny\color{gray},
  numbersep=8pt,
  frame=single,
  framesep=6pt,
  breaklines=true,
  breakatwhitespace=true,
  columns=fullflexible,
  keepspaces=true,
  literate=
    {—}{{---}}3
    {–}{{--}}3
    {∧}{{$\land$}}3
    {∨}{{$\lor$}}3
    {→}{{$\rightarrow$}}3
    {∀}{{$\forall$}}3
    {∃}{{$\exists$}}3
    {∈}{{$\in$}}3
    {≠}{{$\neq$}}3
    {≤}{{$\leq$}}3
    {⟨}{{$\langle$}}3
    {⟩}{{$\rangle$}}3
    {⊢}{{$\vdash$}}3
    {³}{{\textsuperscript{3}}}2
    {♄}{{[Saturn]}}3
    {♃}{{[Jupiter]}}3
    {♂}{{[Mars]}}3
    {♀}{{[Venus]}}3
    {☿}{{[Mercury]}}3
}

\newcommand{\Monas}{\emph{Monas Hieroglyphica}}
\newcommand{\lean}{\textsc{Lean}}

\title{\bfseries Formalising John Dee's \Monas{}\\ in \lean{}~4}
\author{Ivan Matychyn\\[2pt]
	\small Faculty of Mathematics and Computer Science\\
	\small University of Warmia and Mazury in Olsztyn, Poland\\
	\small\texttt{ivan.matychyn@uwm.edu.pl}}
\date{}

\begin{document}
\maketitle

\begin{abstract}
\noindent
John Dee's \Monas{} (Antwerp, 1564) contains twenty-four numbered theorems
about a single composite alchemical glyph, presented in the manner of Euclid:
theorem, demonstration, corollary. In this note the following question is
considered: which of Dee's theorems can be stated precisely in a proof
assistant, and which of them can then be proved? Working in \lean{}~4
(v4.35.0) with the core library only (without \textsc{Mathlib}), we divide the
treatise into three layers: the \emph{syntax} of the glyph, encoded as
inductive data; Dee's \emph{arithmetic}, verified by kernel computation; and
the \emph{geometry} of the final theorems, in connection with which we prove a
classical result not stated by Dee: no ratio of whole numbers doubles the
cube, that is, the cube root of two is irrational (the Delian obstruction).
The accompanying file \texttt{Monas.lean} compiles with no errors and no
\texttt{sorry}. Two observations are made. The cabalistic computation of
Theorem~XVII, according to which the Cross, ``by the best kabbalistic
computation,'' signifies the number 252, is verified, and its proof depends
on no axioms: it is a pure computation. The only theorem of the formalisation
with genuine mathematical content, the Delian obstruction, is proved
constructively and depends only on the axioms \texttt{propext} and
\texttt{Quot.sound}.
\end{abstract}

\tableofcontents

\section{Introduction}

In January 1564, in Antwerp, the mathematician, astronomer and magus John Dee
(1527--1608/9) completed in twelve days a short Latin treatise that he
considered, for the rest of his life, among his most important works
\cite{Josten1964,Peterson}. The \Monas{} (``The Hieroglyphic Monad'') is
devoted to a single sign: a composite glyph built from the semicircle of the
Moon, the circle of the Sun, the cross of the four Elements, and the sign of
Aries. Around this one symbol Dee formulated twenty-four numbered \emph{theorems},
each with a ``demonstration,'' in open imitation of the geometrical method:
``It is by the straight line and the circle that the first and most simple
example and representation of all things may be demonstrated'' (Theorem~I).
Dee claimed that the work would bring a complete change to astronomy,
alchemy, mathematics, linguistics, music, optics and magic; its interpretation
is still disputed \cite{Josten1964,LiteratureAtlas}.

\begin{figure}[t]
\centering
\begin{tikzpicture}
  % the Monad, from top to bottom: Moon, Sun, Cross, Aries
  \begin{scope}[scale=1.35]
    \draw[glyph] (-1.0,0.4) arc (180:360:1.0);
    \draw[glyph] (0,0.05) circle (0.42);
    \fill (0,0.05) circle (0.05);
    \draw[glyph] (0,-0.37) -- (0,-2.6);
    \draw[glyph] (-0.55,-1.3) -- (0.55,-1.3);
    \draw[glyph] (0,-2.6) arc (0:180:0.28) (0,-2.6) arc (180:0:0.28);
  \end{scope}
  % labels
  \begin{scope}[font=\small, anchor=west]
    \draw[dotted] (1.45,0.54) -- (3.5,1.0);
    \node at (3.6,1.0) {\texttt{crescent}: the Moon};
    \draw[dotted] (0.6,0.07) -- (3.5,0.0);
    \node at (3.6,0.0) {\texttt{circle}: the Sun, with its centre};
    \draw[dotted] (0.78,-1.76) -- (3.5,-1.76);
    \node at (3.6,-1.76) {\texttt{line}: the Cross, the sign of the Elements};
    \draw[dotted] (0.8,-3.5) -- (3.5,-3.5);
    \node at (3.6,-3.5) {\texttt{aries}: two crescents with one common point};
  \end{scope}
\end{tikzpicture}
\caption{The Monad, schematic. The list \texttt{monas} of the formalisation is
\texttt{[crescent, circle, line, aries]}, read from top to bottom.}
\label{fig:monas}
\end{figure}
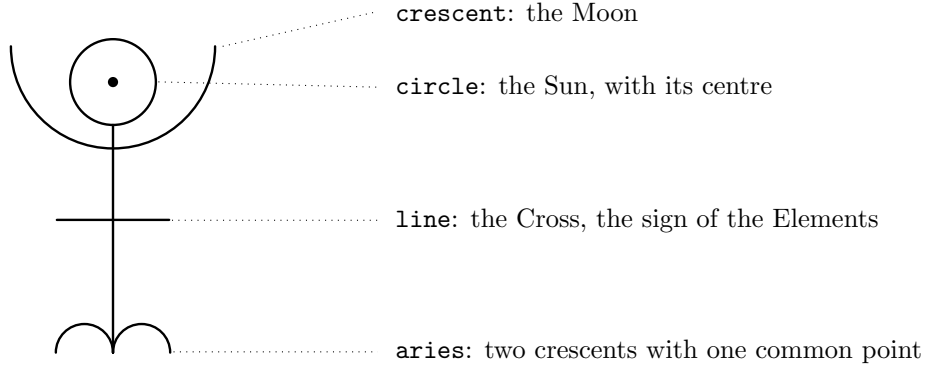

Figure~\ref{fig:monas} shows the glyph and the parts that the formalisation
distinguishes. The drawings in this note are schematic: they were made for this
text from the definitions of Section~\ref{sec:syntax}. They are not facsimiles of
Dee's woodcut, and the proportions are only indicative.

The present note lies between the history of ideas and formal
verification. Dee presents his statements in the form of theorems with
proofs. A proof assistant makes it possible to test this form: a statement
either admits a formalisation, and is then proved or refuted, or it does not
admit one. We therefore consider the theorems one by one and ask whether each
of them can be stated in \lean{}~4 and, if so, whether it holds.

It is shown that the theorems fall into three layers
(Section~\ref{sec:method}): (i) stipulations about the composition of a
symbol, (ii) arithmetic identities, all of which are correct, and
(iii) geometric proportions, for which Dee states no mathematical claim. For
the third layer we prove a related classical result, which is not contained
in the treatise. This division may be of some interest for the commentary
literature on the \Monas{}.

The complete formalisation accompanies this note as the file
\texttt{Monas.lean}. It uses only the \lean{}~4 core library, so it can be
checked by the following command, without need to download a multi-gigabyte \textsc{Mathlib}:
\begin{center}
  \texttt{lean Monas.lean}
\end{center}
Verified with \lean{}~v4.35.0; the output of the final axiom audit is
reproduced in Section~\ref{sec:results}.

\section{Method: a stratification of the twenty-four theorems}
\label{sec:method}

The \Monas{} combines at least five subjects: geometry, arithmetic,
astronomy, alchemy and cabala. Only the first two admit formal verification.
Within them we distinguish three levels, which correspond to three groups of
theorems.

\begin{enumerate}
\item \textbf{Syntax (Theorems I--XV).} The glyph itself as an object of
  study: which signs are composed of which. These claims become inductive
  data declarations. Once the data is laid down, some of Dee's theorems become
  \emph{definitional} truths (proved by \texttt{rfl}), others become
  membership claims about finite lists (proved by computation).
\item \textbf{Arithmetic (Theorems VI, VIII, XVI, XVII, XX).} Dee's
  number-mystical computations: the tetraktys, the measures of the Cross, the
  number 252, are stated precisely enough to be checked. Every identity he
  asserts in these theorems is verified by the \lean{} kernel, mostly by
  \texttt{decide}.
\item \textbf{Geometry (Theorems XXII--XXIV).} The treatise closes with
  ``the proportions already observed by us in the hieroglyphic construction
  of our Monad, presented in diagrammatic form,'' from which the sign of
  Taurus is finally derived. Dee does not state a mathematical claim here;
  he only exhibits proportions. One reading \cite{Wilson} connects
  Theorem~XX, together with the proportions of Theorem~XXIII, with the
  duplication of the cube (the Delian problem). Following this reading, we
  prove the classical result that no ratio of whole numbers doubles the cube.
  This result is not stated by Dee.
\end{enumerate}

Table~\ref{tab:audit} gives the status of all twenty-four theorems.

\begin{table}[t]
\centering\small
\begin{tabular}{@{}lll@{}}
\toprule
Theorem(s) & Content & Status in \texttt{Monas.lean}\\
\midrule
I--V & point, line, circle, Moon, first day & syntax (definitions)\\
VI & Ternary, Septenary, Octad & verified arithmetic\\
VII & Elements resolved into lines & motivates \texttt{atoms}\\
VIII & tetraktys; $X$ is the 21st letter & verified arithmetic\\
IX & lines, circles and the Monas & not formalised (geometric)\\
X & the Monas contains Mercury & verified (computation)\\
XI & Aries is two semicircles & definitional (\texttt{rfl})\\
XII & planetary signs from Moon, Sun, Elements, Aries & verified (computation)\\
XIII & anatomy of Mars, Venus, Mercury & data (the \texttt{parts} table)\\
XIV & the magistery depends on Sun and Moon & not formalised (alchemical)\\
XV & Taurus shows alpha twice & definitional (\texttt{rfl})\\
XVI & quinary, circular numbers, 100, 2500 & verified arithmetic\\
XVII & the Cross signifies 252 & verified arithmetic (no axioms)\\
XVIII--XIX & the egg of the philosophers & not formalised (alchemical)\\
XX & Binary $+$ point $=$ Ternary & verified arithmetic\\
XXI & the Monad redistributed & verified (computation)\\
XXII--XXIV & proportions of the Cross; Taurus & related result proved\\
\bottomrule
\end{tabular}
\caption{Audit of the twenty-four theorems. The connection of the Delian
obstruction with the treatise follows the reading of~\cite{Wilson}, which
relates the duplication of the cube to Theorems~XX and~XXIII.}
\label{tab:audit}
\end{table}

\section{Layer I: the glyph as data}
\label{sec:syntax}

The primitive marks of the hieroglyphic art become an inductive type;
Theorem~XI's ``two marks meeting in one common point'' becomes a binary
constructor:

\begin{lstlisting}[style=lean,numbers=none]
inductive Mark where
  | point : Mark    -- the Monad: the point from which all things emerge
  | line : Mark     -- the straight line; crossed, the sign of the Elements
  | circle : Mark   -- the full circle: the Sun, "a circle having a visible center"
  | crescent : Mark -- the semicircle: the Moon, "the figure of the Horns"
  | binary : Mark -> Mark -> Mark -- two marks meeting in one common point
deriving DecidableEq, Repr
\end{lstlisting}

\begin{figure}[t]
\centering
\begin{tikzpicture}
  \pic[scale=1.5] at (0,0)    {mpoint};
  \pic[scale=1.5] at (2.6,0)  {mline};
  \pic[scale=1.5] at (5.2,0)  {mcircle};
  \pic[scale=1.5] at (7.8,0)  {mcrescent};
  \pic[scale=1.5] at (10.4,0) {maries};
  \pic[scale=1.5] at (13.0,0) {mtaurus};
  \foreach \x/\n/\d in {0/point/the Monad, 2.6/line/the Cross,
      5.2/circle/the Sun, 7.8/crescent/the Moon,
      10.4/aries/two crescents, 13.0/taurus/circle and crescent} {
    \node[font=\small\ttfamily] at (\x,-1.2) {\n};
    \node[font=\footnotesize] at (\x,-1.7) {\d};
  }
  \draw[decorate,decoration={brace,mirror,amplitude=4pt}] (-0.6,-2.2) -- (8.5,-2.2);
  \draw[decorate,decoration={brace,mirror,amplitude=4pt}] (9.3,-2.2) -- (14.1,-2.2);
  \node[font=\footnotesize] at (3.95,-2.85) {the constructors without arguments};
  \node[font=\footnotesize] at (11.7,-2.85) {two uses of \texttt{binary}};
\end{tikzpicture}
\caption{The marks of the type \texttt{Mark}. In the Lean file the mark \texttt{line}
is the straight line; as a part of a sign it appears crossed (the Cross).}
\label{fig:marks}
\end{figure}
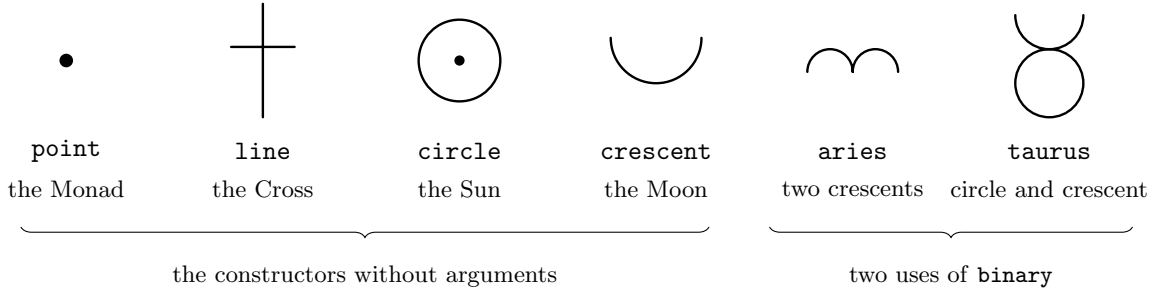

Figure~\ref{fig:marks} draws the four primitive marks and the two signs that
the Lean file builds from them with the constructor \texttt{binary}.

We note one decision of the model: the sign of Aries is \emph{not}
a new primitive but a definition, exactly as Theorem~XI asserts
(``composed of two semicircles connected by one common point''):

\begin{lstlisting}[style=lean,numbers=none]
def aries : Mark := Mark.binary Mark.crescent Mark.crescent
\end{lstlisting}

The function \texttt{atoms} resolves every mark ``into elementary form''
(Theorem~VII's own phrase): composite marks are flattened to the primitive
alphabet. The five planetary characters are then a table, and Theorem~XII's
compositional claim: that the ancient Magi composed all five planetary signs
out of the signs of the Moon, the Sun, the Elements and Aries, is a finite
membership statement:

\begin{lstlisting}[style=lean,numbers=none]
theorem theoremXII :
    forall p : Planet, forall m in sigil p,
      m = Mark.crescent or m = Mark.circle or m = Mark.line
\end{lstlisting}

(In this excerpt the symbols $\forall$, $\in$ and $\lor$ of the \lean{}
source are written in words.) The proof is a case analysis over all marks and all planets,
discharged by \texttt{simp} over the tables; \lean{} reduces each goal to
\texttt{True} or to a contradiction. Two further containment claims of the
treatise are proved in the same way: Theorem~X (the complete sign of Mercury
is present within the Monas) and Theorem~XXI (every planetary character can
be recovered from the Monas; in the formalisation this means that every
elementary mark of every planet occurs among the elementary marks of the
Monas).

\begin{figure}[t]
\centering
\begin{tikzpicture}
	% signs (standard symbols from the package wasysym)
	\foreach \x/\sym in {0/\saturn, 3.1/\jupiter, 6.2/\mars, 9.3/\venus, 12.4/\mercury}
	\node[scale=3.2] at (\x,0) {\sym};
	% names
	\foreach \x/\n in {0/Saturn, 3.1/Jupiter, 6.2/Mars, 9.3/Venus, 12.4/Mercury}
	\node[font=\small] at (\x,-1.2) {\n};
	% entries of the table `parts'
	\pic[scale=0.55] at (-0.45,-2.1) {mcrescent};
	\node[font=\small] at (0,-2.1) {$+$};
	\pic[scale=0.55] at (0.45,-2.1) {mline};
	\pic[scale=0.55] at (2.65,-2.1) {mcrescent};
	\node[font=\small] at (3.1,-2.1) {$+$};
	\pic[scale=0.55] at (3.55,-2.1) {mline};
	\pic[scale=0.55] at (5.75,-2.1) {mcircle};
	\node[font=\small] at (6.2,-2.1) {$+$};
	\pic[scale=0.55] at (6.65,-2.1) {maries};
	\pic[scale=0.55] at (8.85,-2.1) {mcircle};
	\node[font=\small] at (9.3,-2.1) {$+$};
	\pic[scale=0.55] at (9.75,-2.1) {mline};
	\pic[scale=0.55] at (11.55,-2.1) {mcircle};
	\node[font=\small] at (11.98,-2.1) {$+$};
	\pic[scale=0.55] at (12.4,-2.1) {mcrescent};
	\node[font=\small] at (12.82,-2.1) {$+$};
	\pic[scale=0.55] at (13.25,-2.1) {mline};
\end{tikzpicture}
\caption{The five planetary signs (upper row, in their usual shapes) and their
	entries in \texttt{parts} (lower row). The marks are those of
	Figure~\ref{fig:marks}.}
\label{fig:planets}
\end{figure}
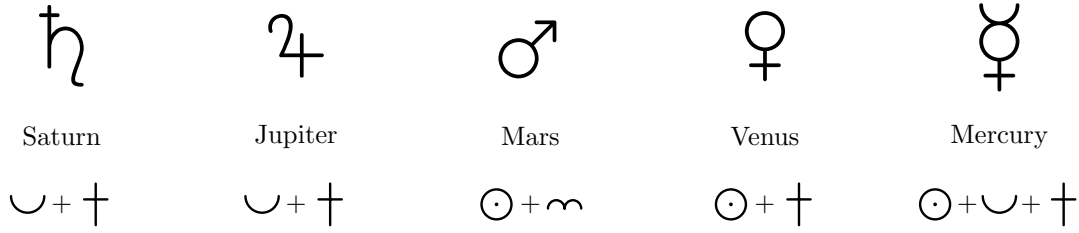

Figure~\ref{fig:planets} shows the five planetary signs, each above its entry
in the table \texttt{parts}. The table does not record the position of the
marks. Hence Saturn and Jupiter have the same entry, and the formalisation
cannot distinguish these two signs.

We make the following remark on method.

\begin{remark}
Once the syntax is laid down, Theorems~XI and~XV are proved by
\texttt{rfl}: they are true \emph{by construction}. The same holds, in a
weaker form, for Theorems~X, XII and~XXI: they are proved by computation, but
they are true because of the way the tables \texttt{parts} and \texttt{monas}
were written. This is not a defect of
the formalisation. It shows that these theorems are stipulations about a
notation and not statements about a subject matter. A proof assistant
distinguishes statements that are true by definition, true by computation,
and true by a non-trivial argument. This distinction is useful for the
reading of a treatise in which the three kinds of statement are not
separated.
\end{remark}

\section{Layer II: the arithmetic of the Cross}
\label{sec:arithmetic}

Dee's numerology is based on two correspondences: the position of letters in
the Latin alphabet as a humanist of 1564 wrote it (no J, U or W; $X$ is the
twenty-first letter, $V$ the twentieth and the fifth of the vowels), and the
Roman numeral values of the letters hidden in the Cross ($V=5$, $L=50$,
$X=10$, $C=100$). Both are encoded as functions, and every arithmetic identity
asserted by Dee is then proved by the tactic \texttt{decide}. Theorem~VI's Ternary, Septenary
and Octad; Theorem~VIII's Pythagorean $1+2+3+4=10$; Theorem~XX's ``the point
must of necessity be present, because with the Binary it constitutes our
Ternary'' ($2+1=3$, and the Quaternary concealed within: $3+1=4$); all of
Theorem~XVI, including the number-theoretic remark that 5 and 25
are \emph{circular numbers}: their squares end in themselves ($25$ ends in
$5$; $625$ ends in $25$), verified as $25 \bmod 10 = 5$ and
$625 \bmod 100 = 25$.

The main statement of this layer is Theorem~XVII, quoted here in full from the
Hamilton-Jones translation \cite{HamiltonJones}:

\begin{quote}\itshape
We may safely conclude, by the best kabbalistic computation, that our Cross,
by a marvellous metamorphosis, may signify for the Initiates two hundred and
fifty-two. Thus: four times five, four times fifty, ten, twenty-one and one,
which added together make two hundred and fifty-two.
\end{quote}

The four right angles of the quinary, four times the letter of fifty, the
Decad, $X$ as the twenty-first letter, and the Monad:

\begin{lstlisting}[style=lean,numbers=none]
theorem theoremXVII :
    4 * romanValue 'V' + 4 * romanValue 'L' + romanValue 'X'
      + latinPos 'X' + 1 = 252 := by decide
\end{lstlisting}

Dee's computation is correct. We note that \texttt{theoremXVII} does not
depend on any axioms (Section~\ref{sec:results}). Formally, the statement is
an equality of two closed terms of type \texttt{Nat}, which the kernel checks
by evaluation.

\section{Layer III: the proportions of the Cross}
\label{sec:delian}

The final three theorems present ``in diagrammatic form the proportions''
of the Cross, in which the parts are made ``equal and unequal,'' and from
which the hieroglyphic figure of Taurus emerges. The interpretation of these
measures is the most disputed part of the treatise. One reading
\cite{Wilson} connects Theorem~XX, together with the proportions of
Theorem~XXIII, with the duplication of the cube, a well-known problem of the
geometry of Dee's time, whose exact solution requires the ratio
$\sqrt[3]{2}$.

Whether or not Dee intended this connection, it leads to a precise statement
that can be proved. We formalise the following classical result, which is not
stated by Dee:

\begin{lstlisting}[style=lean,numbers=none]
theorem impossible :
    forall p q : Nat, q <> 0 -> p * p * p <> 2 * (q * q * q)
\end{lstlisting}

No rational proportion doubles the cube; the cube root of two is irrational.
The proof in \texttt{Monas.lean} is the classical infinite descent, organised
as strong induction on a bound for $p$. From $p^3 = 2q^3$ the root $p$ must be
even, $p = 2r$ (an odd cube is odd, and cannot equal twice anything); then
$q^3 = 4r^3$ forces $q = 2s$ in the same way; and $r^3 = 2s^3$ is a strictly
smaller solution, which is impossible. We note two details of the
implementation, which show what can be done with the \lean{} core library
alone:

\begin{itemize}
\item \emph{Parity without \textsc{Mathlib}.} ``An odd cube is odd'' is proved
  as follows. The cube $(2t+1)^3$ is expanded by the distributivity lemmas of
  the core library, and the result is passed to the tactic \texttt{omega}.
  This tactic solves only linear problems; it treats the products
  $t\cdot t\cdot t$ and $t\cdot t$ as new unknowns $a$ and $b$, and finds that
  the equation $8a + 12b + 6t + 1 = 2Y$, where $Y$ is a natural number, has no
  solution, since its left side is odd and its right side is even.
\item \emph{Normalisation by \texttt{ac\_rfl}.} The identity
  $(2k)(2k)(2k) = 2\cdot2\cdot2\cdot(kkk)$ is proved by \texttt{ac\_rfl},
  which uses only the associativity and commutativity of multiplication.
  After this step the products $q\cdot q\cdot q$, $r\cdot r\cdot r$ and
  $s\cdot s\cdot s$ occur in the equations of the descent only as whole terms.
  The tactic \texttt{omega} treats each of them as a new unknown, so these
  equations become linear, and \texttt{omega} solves them.
\end{itemize}

The Lean file also contains the theorem \texttt{bounds}, which verifies by cross
multiplication the classical bounds $5/4 < \sqrt[3]{2} < 4/3$. That every
ratio of whole numbers is either smaller or greater than $\sqrt[3]{2}$ is the
content of the theorem \texttt{impossible}, not of these bounds.

\section{Results}
\label{sec:results}

The file \texttt{Monas.lean} (324 lines, including the commentary) compiles
with \lean{}~v4.35.0 core, with no errors, no warnings, and no
\texttt{sorry}. The final axiom audit prints:

\begin{lstlisting}[style=lean,numbers=none]
'Monas.theoremXVII' does not depend on any axioms
'Monas.Delian.impossible' depends on axioms: [propext, Quot.sound]
\end{lstlisting}

\noindent We make two remarks.

\begin{remark}
The passage of the treatise with the strongest mystical content is, in the
technical sense, the most trivial one: its verification is a computation that
uses no axioms. Dee's arithmetic is correct; the cabalistic interpretation of
the numbers lies outside the scope of formal verification.
\end{remark}

\begin{remark}
The proof of \texttt{impossible} is constructive: the axioms
\texttt{propext} and \texttt{Quot.sound} are compatible with intuitionistic
logic, and the law of the excluded middle is not used. An infinite descent of
this kind is a constructive argument, as are the classical Greek proofs of
irrationality.
\end{remark}

\section{Future work}

We indicate three possible continuations.
\begin{enumerate}
\item \emph{Constructive geometry.} The exact gnomon proportions of
  Theorem~XXIII could be formalised as a statement in the real numbers,
  following the classical Greek programme for the Delian problem (two mean
  proportionals). This requires \textsc{Mathlib}, since the core library
  contains no ordered field. Such a formalisation would replace the
  interpretation of Section~\ref{sec:delian} by a reconstruction of Dee's own
  diagram.
\item \emph{Circular numbers.} Dee's ``circular numbers'' are the automorphic
  numbers of modern number theory. Classifying the two-digit automorphic
  numbers ($25$ and $76$) is within reach of the Lean file's toolkit and
  would complete the treatment of Theorem~XVI.
\item \emph{Teaching.} The Lean file was written with a view to a \lean{}~4
  course in preparation (December 2026). On a single example it uses
  definitions as data, \texttt{rfl}, \texttt{decide}, \texttt{omega},
  induction and infinite descent. The fact that a theorem of 1564 depends on
  no axioms may serve as an introduction to the subject in a lecture.
\end{enumerate}

\appendix

\section{The complete \lean{} source}
\label{app:source}

The complete file \texttt{Monas.lean}, as verified with \lean{}~v4.35.0
(core library only). Comments preserve the quotations from the
Hamilton-Jones translation used as proof documentation; the planetary glyphs
in the \texttt{parts} table are transliterated for the typesetter.

\begin{lstlisting}[style=lean]
/-!
# John Dee's *Monas Hieroglyphica* (Antwerp, 1564) — a Lean 4 formalisation

John Dee's *Monas Hieroglyphica* ("The Hieroglyphic Monad"), written in twelve
days and published at Antwerp in 1564, unfolds twenty-four theorems about a
single composite glyph — the *Monas* — built from the signs of the Moon, the
Sun, the four Elements and Aries. Dee presents it in the style of Euclid:
numbered theorems with demonstrations, applied not to abstract figures but to
his hieroglyph.

This file takes that Euclidean posture seriously: which of Dee's theorems can
be *stated* precisely, and which can then be *proved*? Three layers emerge.

* **Syntax** — the glyph as data: inductive types for the primitive marks and
  the planetary characters (Theorems I–XV).
* **Number** — Dee's own arithmetic, verified exactly as he states it:
  the tetraktys, the measures of the Cross, the number 252 (Theorems VI, VIII,
  XVI, XVII).
* **Geometry** — the proportions of the Cross with which the treatise closes
  (Theorems XXII–XXIV). Following the reading that connects these measures
  with the duplication of the cube — the Delian problem — we verify the
  classical obstruction: *no rational proportion doubles the cube*
  (`Delian.impossible`). That last theorem is not stated by Dee; it is a classical result connected
  with his construction under this reading.

Quotations follow J. W. Hamilton-Jones's 1947 English translation; the
scholarly standard is C. H. Josten's translation and commentary, *Ambix* XII
(1964), 84–221. Verified with Lean 4 core alone (no Mathlib required).
-/

namespace Monas

/-! ## Layer 0: the syntax of the glyph -/

/-- The primitive marks of the hieroglyphic art (Theorems I–V). -/
inductive Mark where
  | point : Mark    -- the Monad: the point from which all things emerge (Thms II–III)
  | line : Mark     -- the straight line; crossed, the sign of the Elements (Thms I, VII)
  | circle : Mark   -- the full circle: the Sun, "a circle having a visible center" (Thm III)
  | crescent : Mark -- the semicircle: the Moon, "the figure of the Horns" (Thms IV–V)
  | binary : Mark → Mark → Mark -- two marks meeting in one common point (Thm XI)
deriving DecidableEq, Repr

/-- **Thm XI.** "The mystical sign of the Ram, composed of two semicircles
connected by one common point." -/
def aries : Mark := Mark.binary Mark.crescent Mark.crescent

/-- **Thm XV** (note). "By a simple transposition of place we show the letter
alpha twice, by a circle and a half-circle, being simply tangents which touch
one another." -/
def taurus : Mark := Mark.binary Mark.circle Mark.crescent

/-- Every mark resolved "into elementary form" (Thm VII): the composite marks
decomposed into the primitive alphabet of point, line, circle and crescent. -/
def atoms : Mark → List Mark
  | .point   => [.point]
  | .line    => [.line]
  | .circle  => [.circle]
  | .crescent => [.crescent]
  | .binary a b => atoms a ++ atoms b

/-- The four right angles of the elemental Cross (Thms VI–VII). -/
inductive Element where
  | fire | air | water | earth
deriving DecidableEq, Repr

/-- The wandering stars of the inferior astronomy (Thms XII–XIII). -/
inductive Planet where
  | saturn | jupiter | mars | venus | mercury
deriving DecidableEq, Repr

/-- **Thm XII.** "The very ancient wise men and Magi have transmitted to us
five hieroglyphical signs of the planets, all of which are composed out of the
signs used for the Moon and the Sun, together with the sign of the Elements
and the hieroglyphical sign of Aries." -/
def parts : Planet → List Mark
  | .saturn  => [.crescent, .line]          -- ♄ : the lunar crescent upon the cross
  | .jupiter => [.crescent, .line]          -- ♃ : the cross beneath the crescent
  | .mars    => [.circle, aries]             -- ♂ : "formed from the hieroglyphs of the Sun and Aries" (Thm XIII)
  | .venus   => [.circle, .line]             -- ♀ : "produced from that of the Sun and the Elements" (Thm XIII)
  | .mercury => [.circle, .crescent, .line]  -- ☿ : "by the complete Lunar and Solar magic of the Elements" (Thm XIII)

/-- A planetary character resolved into its elementary marks. -/
def sigil (p : Planet) : List Mark := (parts p).flatMap atoms

/-- The complete Monad, read from top to bottom (Thm X): the crescent of the
Moon upon the circle of the Sun upon the elemental Cross, with the sign of
Aries at its base. -/
def monas : List Mark := [Mark.crescent, Mark.circle, Mark.line, aries]

/-- The Monas resolved into its elementary marks. -/
def monasSigil : List Mark := monas.flatMap atoms

/-! ## Layer 1: structural theorems -/

/-- **Thm XI.** "The mystical sign of the Ram, composed of two semicircles
connected by one common point": resolved into elementary form, the sign of
Aries is exactly two semicircles. -/
theorem theoremXI : atoms aries = [Mark.crescent, Mark.crescent] := rfl

/-- **Thm XV** (note). The figure of Taurus shows alpha twice: a circle and a
half-circle, tangents touching at one common point. -/
theorem theoremXV : atoms taurus = [Mark.circle, Mark.crescent] := rfl

/-- **Thm XII.** Every planetary character, resolved into elementary form, is
composed only of the marks of the Moon, the Sun and the Elements — Dee's
compositional claim about the five characters of the planets, verified by
computation. -/
theorem theoremXII :
    ∀ p : Planet, ∀ m ∈ sigil p,
      m = Mark.crescent ∨ m = Mark.circle ∨ m = Mark.line := by
  intro p m hm
  cases m <;> cases p <;> simp [sigil, parts, atoms, aries] at hm ⊢

/-- **Thm X.** "It is the Mercury of the Philosophers, the greatly celebrated
microcosm and ADAM": every elementary mark of Mercury is present within the
Monas. -/
theorem theoremX : ∀ m ∈ sigil Planet.mercury, m ∈ monasSigil := by
  intro m hm
  cases m <;> simp [sigil, parts, atoms, aries, monas, monasSigil] at hm ⊢

/-- **Thm XXI.** "We distribute our Monad, now looked at from a different
aspect, and analysed in a different manner": every planetary character can be
recovered from the Monas. -/
theorem theoremXXI : ∀ p : Planet, ∀ m ∈ sigil p, m ∈ monasSigil := by
  intro p m hm
  cases m <;> cases p <;> simp [sigil, parts, atoms, aries, monas, monasSigil] at hm ⊢

/-! ## Layer 2: the arithmetic of the Cross -/

/-- Dee's Latin alphabet, as a humanist of 1564 writes it: no J, U or W.
X is the twenty-first letter, V the twentieth. -/
def latinPos (c : Char) : Nat := match c with
  | 'A' => 1 | 'B' => 2 | 'C' => 3 | 'D' => 4 | 'E' => 5
  | 'F' => 6 | 'G' => 7 | 'H' => 8 | 'I' => 9 | 'K' => 10
  | 'L' => 11 | 'M' => 12 | 'N' => 13 | 'O' => 14 | 'P' => 15
  | 'Q' => 16 | 'R' => 17 | 'S' => 18 | 'T' => 19 | 'V' => 20
  | 'X' => 21
  | _ => 0

/-- The Roman numeral values of the letters hidden in the Cross (Thm XVI):
"the fifth of the vowels" for five, the letter of fifty, the Decad for X. -/
def romanValue (c : Char) : Nat := match c with
  | 'V' => 5
  | 'L' => 50
  | 'X' => 10
  | 'C' => 100
  | _ => 0

/-- **Thm VI.** The Ternary is "two straight lines having a copulative
center"; the Quaternary is "produced by the four straight lines"; together
they manifest the Septenary, and "either of these elements, the lines or the
right angles, repeated twice, afford us in the most secret manner the Octad". -/
theorem theoremVI :
    2 + 1 = 3          -- two lines and their copulative point: the Ternary
    ∧ 3 + 4 = 7        -- Ternary and Quaternary: the first manifested Septenary
    ∧ 2 * 4 = 8        -- the lines or the right angles, repeated twice: the Octad
    := by decide

/-- **Thm VIII.** "Pythagoras was in the habit of saying: 1 + 2 + 3 + 4 make
10. It is not by chance that the right-angled Cross — the twenty-first letter
of the Roman alphabet — was taken by the most ancient of the Roman
philosophers to represent the Decad." -/
theorem theoremVIII :
    1 + 2 + 3 + 4 = 10 ∧ latinPos 'X' = 21 := by decide

/-- **Thm XX.** "The point must of necessity be present, because with the
Binary it constitutes our Ternary" — and the Quaternary is concealed within
the Ternary. -/
theorem theoremXX : 2 + 1 = 3 ∧ 3 + 1 = 4 := by decide

/-- **Thm XVI.** The quinary is "the exact half of our Decad"; from the
twentieth letter of the alphabet and the fifth of the vowels comes 25, the
first of the *circular numbers* — numbers whose square ends in themselves.
The two parts of the Cross produce the number one hundred and, "by the law
of squares", 2500; "we are thereby induced to utilise this progression:
one — ten — one hundred." -/
theorem theoremXVI :
    5 * 2 = 10                                  -- the quinary: half of the Decad
    ∧ latinPos 'V' = 20                         -- V: twentieth letter of the alphabet
    ∧ romanValue 'V' = 5                        -- V: fifth of the vowels
    ∧ latinPos 'V' + romanValue 'V' = 25        -- "the number twenty-five is produced"
    ∧ 5 * 5 = 25 ∧ 25 * 25 = 625                -- squares of the circular numbers
    ∧ 25 % 10 = 5                               -- 5 is circular: its square ends in 5
    ∧ 625 % 100 = 25                            -- 25 is circular: its square ends in 25
    ∧ romanValue 'L' + romanValue 'L' = 100     -- the two parts: the number one hundred
    ∧ romanValue 'L' * romanValue 'L' = 2500    -- "by the law of squares"
    ∧ 10 * 1 = 10 ∧ 10 * 10 = 100               -- one — ten — one hundred
    := by decide

/-- **Thm XVII.** "We may safely conclude, by the best kabbalistic
computation, that our Cross, by a marvellous metamorphosis, may signify for
the Initiates two hundred and fifty-two. Thus: four times five, four times
fifty, ten, twenty-one and one, which added together make two hundred and
fifty-two."

(four right angles of the quinary; four times the letter of fifty; the Decad;
X the twenty-first letter; the Monad) — Dee's computation, machine-checked. -/
theorem theoremXVII :
    4 * romanValue 'V' + 4 * romanValue 'L' + romanValue 'X' + latinPos 'X' + 1
      = 252 := by decide

/-! ## Layer 3: the proportions of the Cross (Thms XXII–XXIV) -/

/- Dee closes the Monas by presenting "in diagrammatic form the proportions
already observed by us in the hieroglyphic construction of our Monad", where
the parts of the Cross are made "equal and unequal", and from which the figure
of Taurus finally emerges. Commentators have connected these measures with the
most famous problem standing before the geometers of Dee's age: the
duplication of the cube — the Delian problem — whose exact solution requires
the proportion 2^(1/3). We verify the classical obstruction: no proportion of
whole numbers can double the cube. (This result is not stated by Dee.)
-/

namespace Delian

/-- Even naturals. -/
def Even (n : Nat) : Prop := ∃ k, n = 2 * k

/-- Odd naturals. -/
def Odd (n : Nat) : Prop := ∃ k, n = 2 * k + 1

theorem even_or_odd (n : Nat) : Even n ∨ Odd n := by
  induction n with
  | zero => exact Or.inl ⟨0, rfl⟩
  | succ n ih =>
    cases ih with
    | inl h =>
      cases h with
      | intro k hk => exact Or.inr ⟨k, by omega⟩
    | inr h =>
      cases h with
      | intro k hk => exact Or.inl ⟨k + 1, by omega⟩

/-- An odd cube is odd: `(2t+1)³` can never equal twice anything. The proof
unfolds the cube by distributivity and hands the linear remainders to
`omega`, which finds the parity contradiction. -/
private theorem odd_cube_ne_two_mul {t Y : Nat} :
    (2 * t + 1) * (2 * t + 1) * (2 * t + 1) ≠ 2 * Y := by
  intro h
  simp only [Nat.mul_add, Nat.mul_one, Nat.one_mul, Nat.two_mul,
    Nat.mul_comm, Nat.mul_left_comm] at h
  omega

/-- The cube of an even number, normalised: `(2k)³ = 8k³`, proved up to
associativity and commutativity of multiplication. -/
private theorem cube_two_mul (k : Nat) :
    (2 * k) * (2 * k) * (2 * k) = 2 * 2 * 2 * (k * k * k) := by ac_rfl

private theorem cube_pos {n : Nat} (h : 0 < n) : 0 < n * n * n :=
  Nat.mul_pos (Nat.mul_pos h h) h

/-- **The Delian obstruction.** There are no naturals `p q` with `q ≠ 0` and
`p³ = 2q³`: the cube cannot be doubled by any rational proportion, and the
cube root of two is irrational.

The proof is the classical infinite descent, organised as strong induction on
a bound for `p`: from `p³ = 2q³` the root `p` must be even, `p = 2r`; then
`q³ = 4r³` forces `q = 2s` in the same way; and `r³ = 2s³` is a strictly
smaller solution — impossible. -/
theorem impossible : ∀ p q : Nat, q ≠ 0 → p * p * p ≠ 2 * (q * q * q) := by
  have aux : ∀ N p q : Nat, p ≤ N → q ≠ 0 → p * p * p = 2 * (q * q * q) → False := by
    intro N
    induction N with
    | zero =>
      intro p q hp hq h
      have hpz : p = 0 := by omega
      rw [hpz] at h
      simp only [Nat.zero_mul] at h
      cases q with
      | zero => exact hq rfl
      | succ k =>
        have hpos : 0 < (k + 1) * (k + 1) * (k + 1) := cube_pos (by omega)
        omega
    | succ N ihN =>
      intro p q hp hq h
      have hpz : p ≠ 0 := by
        intro hpz
        rw [hpz] at h
        simp only [Nat.zero_mul] at h
        cases q with
        | zero => exact hq rfl
        | succ k =>
          have hpos : 0 < (k + 1) * (k + 1) * (k + 1) := cube_pos (by omega)
          omega
      cases even_or_odd p with
      | inr hodd =>
        cases hodd with
        | intro t ht =>
          rw [ht] at h
          exact odd_cube_ne_two_mul h
      | inl heven =>
        cases heven with
        | intro r hr =>
          rw [hr, cube_two_mul] at h
          have hq4 : q * q * q = 2 * (2 * (r * r * r)) := by omega
          cases even_or_odd q with
          | inr hoddq =>
            cases hoddq with
            | intro t ht =>
              rw [ht] at hq4
              exact odd_cube_ne_two_mul hq4
          | inl hevenq =>
            cases hevenq with
            | intro s hs =>
              rw [hs, cube_two_mul] at hq4
              have hfinal : r * r * r = 2 * (s * s * s) := by omega
              exact ihN r s (by omega) (by omega) hfinal
  intro p q hq h
  exact aux p p q (Nat.le_refl p) hq h

/-- Dee's proportions of the Cross are "equal and unequal" indeed: any whole
number ratio falls short of or overshoots the cube root of two. Between 5/4
and 4/3 it lies — and neither of them is it, as one verifies by cross
multiplication. -/
theorem bounds :
    5 * 5 * 5 < 2 * (4 * 4 * 4) ∧ 2 * (3 * 3 * 3) < 4 * 4 * 4 := by decide

end Delian

end Monas

#print axioms Monas.theoremXVII
#print axioms Monas.Delian.impossible
\end{lstlisting}

\section{Reproduction}
\label{app:repro}

The Lean file is self-contained and requires no \textsc{Mathlib}. With any
recent \lean{}~4 toolchain:
\begin{center}
  \texttt{lean Monas.lean}
\end{center}
The build is silent on success except for the two \texttt{\#print axioms}
audits reproduced in Section~\ref{sec:results}. Verified toolchain: \lean{}
v4.35.0, commit \texttt{eb6cafb} (official release binaries). No options or
build system are needed.

\end{document}